# Forced Symmetry-Breaking via Boundary Conditions

Zhen Mei *

Fachbereich Mathematik

Universität Marburg

D-35032 Marburg, Germany

Chih-Wen Shih†

Department of Applied Mathematics

National Chiao-Tung University

Hsinchu, Taiwan, R.O.C.

March 9, 1998

## Abstract

We study impact of a forced symmetry-breaking in boundary conditions on the bifurcation scenario of a semilinear elliptic partial differential equation. We show that for the square domain the orthogonality of eigenfunctions of the Laplacian may compensate partially the loss of symmetries in the boundary conditions and allows some solution to have more symmetries than the imposed boundary conditions.

**AMS(MOS) Subject Classification:** 58G28, 65J15, 35B32

---

*Partially supported by the DFG of Germany via grant Me 1336/1-3 and by the Fields Institute for Research in Mathematical Sciences, Canada.

†Work supported by the National Science Council of R.O.C.

# 1 Introduction

Boundary conditions have sophisticate influence on behavior of reaction-diffusion equations. Since diffusion is the underlying mechanism for the spatial pattern formation in chemical reactions, spatial structure of solutions of reaction-diffusion equations can be sensitive to boundary conditions. Posing and realizing appropriate boundary conditions, e.g. for chemical reactions in open and large systems, is delicate (cf. Gray/Scott [7]). Furthermore, components of a system of reaction-diffusion equations can be imposed with different boundary conditions. These make a big difference in behavior of systems and scalar equations, e.g. in Hopf bifurcations, spirals and other pattern formations. Typically multiple bifurcations occur more likely in systems. Moreover, stability of the bifurcating solution branches varies considerably from a scalar equation to a system. To distinguish influence of boundary conditions from that of interactions among the different species (components) of systems, we consider a scalar reaction-diffusion equation

$$\Delta u + \lambda u = f(u, \lambda) \quad \text{in} \quad \Omega := (0, \pi) \times (0, \pi). \tag{1}$$

We assume that the mapping $f : \mathbf{R} \times \mathbf{R} \to \mathbf{R}$ is sufficiently smooth and

$$f(0, \lambda) = 0, \qquad D_u f(0, \lambda) = 0, \tag{2}$$

i.e., it describes the nonlinearity of the problem (1) and implies that

$$u \equiv 0, \quad \lambda \in \mathbf{R}$$

is a trivial solution of (1).

We are interested in impact of symmetry-breaking in boundary conditions on bifurcation scenarios. To this end we consider a square domain and impose the following conditions along its four sides

$$\begin{aligned}
h_0(\mu)u(x,0) - h_1(\mu)\frac{\partial u}{\partial y}(x,0) &= 0, \\
h_0(\mu)u(x,\pi) + h_1(\mu)\frac{\partial u}{\partial y}(x,\pi) &= 0, \\
\frac{\partial u}{\partial x}(0,y) &= 0, \\
\frac{\partial u}{\partial x}(\pi,y) &= 0.
\end{aligned} \tag{3}$$

Here $h_0$, $h_1 : [0,1] \to R$ are smooth functions and satisfy

$$h_0(0) = h_1(1) = 0, \ h_0(\mu) \neq 0 \ \text{for} \ \mu \in (0,1], \ h_1(\mu) \neq 0 \ \text{for} \ \mu \in [0,1). \tag{4}$$

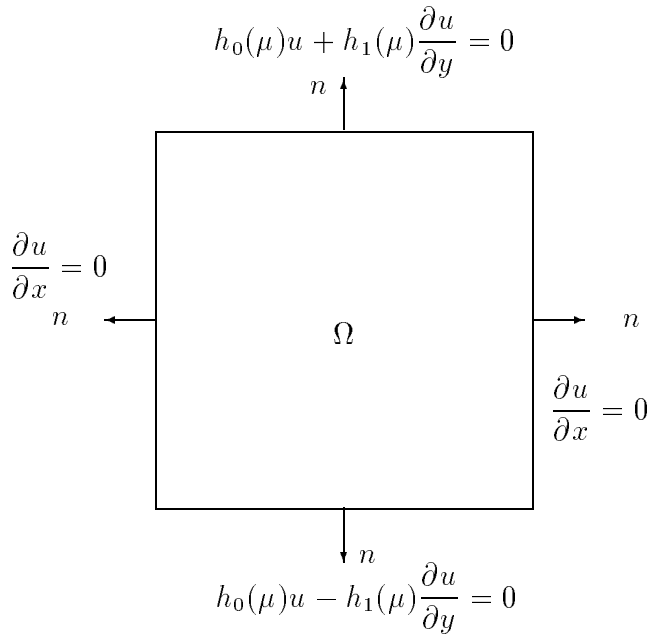

Figure 1: *A homotopy between homogeneous Neumann and mixed boundary conditions.*

The boundary conditions (3) break the $D_4$-symmetry into $D_2$-symmetry. More precisely, properties of $h_0(\mu)$, $h_1(\mu)$ make (3) a homotopy from the homogeneous Neumann boundary conditions along the four sides of $\Omega$ at $\mu = 0$ to the mixed boundary conditions at $\mu = 1$, which are of the Neumann type along the sides $x = 0, \pi$; and of the Dirichlet type at $y = 0, \pi$.

Homotopy of boundary conditions has been used by Fiedler [3] and Gardner [4] to study global attractors and nonsingular solutions of a class of reaction-diffusion equations. They show that these solution sets are independent of boundary conditions. Nevertheless, if the equilibrium is nonhyperbolic and a bifurcation occurs, the bifurcation scenario, e.g. structure of attractors, may vary with respect to boundary conditions. This has been observed by Dillon/Maini/Othmer [2] in the study of pattern formation in generalized one-dimensional Turing systems and by Mei/Theil [13] in the analysis of steady state bifurcations as well as by Holder/Schaeffer [9] and Schaeffer/Golubitsky [14] on mode-jumping of von Kármán equations. Using (3) as an example we study in this paper how reaction-diffusion equations react to a symmetry-breaking in boundary conditions.

An outline of this paper is as follows. In Section 2 we consider variational form and symmetries of the equation (1). Section 3 describes bifurcation points of (1) along the trivial solution curve. In Section 4 the problem (1) at bifurcation points is reduced to algebraic equations via the well-known Liapunov-Schmidt method. We derive the bifurcation scenario at simple and double bifurcation points in Section 5 and illustrate these with a simple example in Section 6.

# 2 Variational Equations and Symmetries

The classical regularity theory of elliptic problems ensures the $C^{2,s}$-Hölder continuity of solutions of the linear problem

$$\Delta u + \lambda u = g$$

with the boundary conditions (3) inside the square $\Omega$ (cf. Wigley [15]). However, differentiability of solutions at the four corners depends strongly on properties of $g$. This linear problem is involved in the analysis of bifurcations of (1) with various right hand sides. here we want to avoid technicalities for the classical solutions and write the problem (1) into variational form. Namely, we study bifurcations of its weak solutions.

## 2.1 Weak form

We consider the Sobolev space

$$X := H^1(\Omega) = \left\{ u \in L^2(\Omega); \quad \frac{\partial u}{\partial x}, \frac{\partial u}{\partial y} \in L^2(\Omega) \right\}$$

with the norm $\|\cdot\|_{1,\Omega}$, and for $u, v \in X$, $\mu \in [0, 1)$ the parameter-dependent bilinear form

$$b_\mu(u,v) := -\int_\Omega (\nabla u \nabla v + uv)\, dxdy - \frac{h_0(\mu)}{h_1(\mu)} \int_0^\pi [u(x,0)v(x,0) + u(x,\pi)v(x,\pi)]\, dx. \quad (5)$$

For $\mu = 1$ we choose the bilinear form

$$b_1(u,v) := -\int_\Omega (\nabla u \nabla v + uv)\, dxdy, \quad (6)$$

defined in the space $\tilde{X} \times \tilde{X}$ and

$$\tilde{X} := \{u \in H^1(\Omega) \,;\, u \text{ satisfies the boundary conditions (3) for } \mu = 1\}$$

The weak form of the linear problem

$$\Delta u - u = g \quad \text{in } \Omega \quad (7)$$

with the boundary conditions (3) is

$$\text{Find } u \in H^1(\Omega), \quad \text{such that } b_\mu(u,v) = (g,v) \quad \text{for all } v \in H^1(\Omega). \quad (8)$$

In particular, weak form of the Neumann problem corresponds to $\mu = 0$.

For domains $\Omega \in C^{0,1}$, typically rectangle and $L$-domains (cf. Hackbusch [8, pp.118]), we have

$$\|u\|_{H^{\frac{1}{2}}(\partial\Omega)} \leq C\|u\|_{1,\Omega} \quad \text{for all} \quad u \in H^1(\Omega),$$

$$\left\|\frac{\partial u}{\partial n}\right\|_{H^{\frac{1}{2}}(\partial\Omega)} \leq C\|u\|_{2,\Omega} \quad \text{for all} \quad u \in H^2(\Omega),$$

where $C > 0$ is a constant. Note that

$$\|\phi\|_{L^2(\partial\Omega)} \leq \|\phi\|_{H^{\frac{1}{2}}(\partial\Omega)} = \inf_{u|_{\partial\Omega}=\phi} \|u\|_{1,\Omega} \quad \text{for all} \quad \phi \in H^{\frac{1}{2}}(\partial\Omega).$$

Thus the bilinear form $b_\mu(\cdot,\cdot)$ is continuous and coercive on $X \times X$. Moreover, if the inequality $h_0(\mu)h_1(\mu) \geq 0$ holds, it is elliptic. Therefore, the problem (8) has a unique solution $u(\mu) \in X$ for every $\mu \in [0,1)$ and all $g \in H^{-1}(\Omega)$. The solution $u(\mu)$ satisfies the boundary conditions (3).

Denote the solution operator of (8) as

$$T(\mu): g \in H^{-1}(\Omega) \longmapsto T(\mu)g = u(\mu) \in H^1(\Omega). \tag{9}$$

The operator $T(\mu)$ is linear and bounded for any fixed $\mu \in [0,1)$. Furthermore, owing to the symmetry of Laplacian it is self-adjoint. This can be seen from the following equality for all $f, g \in H^{-1}(\Omega)$

$$\langle T(\mu)f, g\rangle_{H^1(\Omega) \times H^{-1}(\Omega)}$$
$$= \int_\Omega (T(\mu)f)g \, dxdy$$
$$= \int_\Omega \left(\nabla(T(\mu)f)\nabla(T(\mu)g) + (T(\mu)f)(T(\mu)g)\right) dxdy$$
$$+ \frac{h_0(\mu)}{h_1(\mu)} \int_0^\pi \left[(T(\mu)f)(T(\mu)g)\Big|_{y=0} + (T(\mu)f)(T(\mu)g)\Big|_{y=\pi}\right] dx$$
$$= \int_\Omega (-\Delta + I)(T(\mu)f)(T(\mu)g) \, dxdy + \int_{\partial\Omega} (T(\mu)g)\frac{\partial}{\partial n}T(\mu)f \, ds$$
$$+ \frac{h_0(\mu)}{h_1(\mu)} \int_0^\pi \left[(T(\mu)f)(T(\mu)g)\Big|_{y=0} + (T(\mu)f)(T(\mu)g)\Big|_{y=\pi}\right] dx$$
$$= \int_\Omega f \cdot (T(\mu)g) \, dxdy$$
$$= \langle f, T(\mu)g\rangle_{H^{-1}(\Omega) \times H^1(\Omega)}.$$

More importantly is the fact that together with the Gelfand relation

$$H^1(\Omega) \overset{c}{\hookrightarrow} L^2(\Omega) \hookrightarrow H^{-1}(\Omega)$$

the Riesz-Schauder theory is applicable to the operator $T(\mu)$ for all $\mu \in [0,1)$, so that $T(\mu)$ has the following properties.

- Spectrum of the operator $T(\mu)$ consists of eigenvalues. There are maximally countable eigenvalues and can be ordered as

$$\lambda_1 \geq \lambda_2 \geq \cdots \to 0.$$

The eigenspace associated to each eigenvalue $\lambda_i, i = 1, 2, \ldots$ is finite dimensional;

- For $i = 1, 2, \ldots$ the equation $T(\mu)u - \lambda u = f$ is solvable if and only if $f \perp \text{Ker}(T(\mu) - \lambda I)$.

For $\mu = 1$ we obtain the same conclusions with the bilinear form (6). In the sequel we consider the weak form

$$G(u, \lambda, \mu) := u + (\lambda + 1)T(\mu)u - T(\mu)f(u, \lambda) = 0. \tag{10}$$

The mapping $G : X \times R \to X$ is obviously as smooth as $f$ in $(u, \lambda)$. It is also continuously differentiable in $\mu$ due to the following properties of $T(\mu)$.

**Lemma 1** *(Mei [12]) The operator $T(\mu)$ is continuous and differentiable with respect to $\mu$ in $[0, 1)$. Furthermore, the derivative $u'(\mu) = T'(\mu)g =: v(\mu)$ for all $g \in Y$ is given as the weak solution of*

$$\begin{aligned}
\Delta v - v &= 0 \quad \text{in } \Omega, \\
\frac{h_0(\mu)}{h_1(\mu)} v(x, 0) - \frac{\partial v}{\partial y}(x, 0) &= -\left(\frac{h_0(\mu)}{h_1(\mu)}\right)' u(x, 0), \\
\frac{h_0(\mu)}{h_1(\mu)} v(x, \pi) + \frac{\partial v}{\partial y}(x, \pi) &= -\left(\frac{h_0(\mu)}{h_1(\mu)}\right)' u(x, \pi) \\
\frac{\partial v}{\partial x}(0, y) = 0, \quad \frac{\partial v}{\partial x}(\pi, y) &= 0.
\end{aligned} \tag{11}$$

To calculate $T(\mu)'g$, we denote $u = T(\mu)g$ the weak solution of the equation $\Delta u - u = g$ with the boundary condition (3). Define

$$\hat{v} := \left(\frac{h_0(\mu)}{h_1(\mu)}\right)' \left(-\frac{y^2}{\pi} + y\right) u(x, y). \tag{12}$$

It is easy to verify that $\hat{v}$ satisfies the boundary conditions in (11). Let $v = w + \hat{v}$ and substitute it into (11). We obtain the equation

$$\Delta w - w = -(\Delta \hat{v} - \hat{v})$$

with the boundary conditions (3). Furthermore,

$$-(\Delta \hat{v} - \hat{v}) = -\left(\frac{h_0(\mu)}{h_1(\mu)}\right)' \left[(-\frac{y^2}{\pi} + y)g - \frac{2}{\pi}u + 2(-\frac{2y}{\pi} + 1)\frac{\partial u}{\partial y}\right].$$

Hence, the weak solution $v$ of the equation (11) is

$$\begin{aligned} v &= -T(\mu)(\Delta \hat{v} - \hat{v}) + \hat{v} \\ &= \left(\frac{h_0(\mu)}{h_1(\mu)}\right)' \left\{ T(\mu) \left[\frac{2}{\pi} T(\mu) g + 2(\frac{2y}{\pi} - 1)\frac{\partial}{\partial y}(T(\mu)g) + \left(\frac{y^2}{\pi} - y\right) g \right] \right. \\ &\quad + \left. \left(-\frac{y^2}{\pi} + y\right) T(\mu) g \right\}. \end{aligned} \qquad (13)$$

**Remark:** Note that $T(\mu)$ is self-adjoint, so is its derivative $T'(\mu)$. Based on the formulation (11), one can calculate the higher order derivatives of $u(\mu)$ in a similar manner.

## 2.2  Symmetries

Let $D_4$ be the dihedral group of the unit square $\Omega$ and

$$S_1(x, y) = (1 - x, y), \qquad R(x, y) = (1 - y, x)$$

be its generators. With $Z_2 := \{1, -1\}$, we define $Z_2 \times D_4 = \{\pm \delta; \; \delta \in D_4\}$ and its actions on $Y := L^2(\Omega) \, (\supset X)$ as

$$\gamma u(x, y) = \pm u(\delta^{-1}(x, y)) \qquad \text{for all } \gamma = \pm \delta, \; \delta \in D_4 \text{ and } u \in Y. \qquad (14)$$

The function spaces $X$, $Y$ are obviously $Z_2 \times D_4$-invariant. Similarly, the $L^2$-product is also $Z_2 \times D_4$-invariant. Corresponding to the boundary condition (3) we are particularly interested in the subgroup

$$D_2 := \{S_1, R^2; \; S_1 R^2, I\}. \qquad (15)$$

Let

$$\Gamma := \begin{cases} Z_2 \times D_4 & \text{if } f(u, \lambda) \text{ is an odd function in } u \text{ and } \mu = 1; \\ D_4 & \text{if } f(u, \lambda) \text{ is not odd in } u \text{ and } \mu = 1; \\ Z_2 \times D_2 & \text{if } f(u, \lambda) \text{ is an odd function in } u \text{ and } \mu \neq 1; \\ D_2 & \text{if } f(u, \lambda) \text{ is not odd in } u \text{ and } \mu \neq 1. \end{cases}$$

The $\Gamma$-*equivariance* of the mapping $G$, i.e.,

$$G(\gamma u, \lambda, \mu) = \gamma G(u, \lambda, \mu) \qquad \text{for all} \quad \gamma \in \Gamma, \; u \in X, \; \lambda \in \mathbf{R}.$$

can be verified directly via the generators $\pm S_1$ and $\pm R$, respectively.

## 3  Bifurcation Points

Since $D_\lambda G(0,\lambda,\mu) = D_\mu G(0,\lambda,\mu) \equiv 0$ and $D_u G(0,\lambda,\mu) = I + (\lambda + 1)T(\mu)$ for all $\lambda,\ \mu \in \mathbf{R}$, a bifurcation occurs at a point $(0,\lambda,\mu)$ on the trivial solution manifold $\{(0,\lambda,\mu);\ \lambda \in \mathbf{R},\ \mu \in [0,1]\}$ of (10) if the linearized problem

$$D_u G(0,\lambda,\mu)u = u + (\lambda+1)T(\mu)u = 0 \tag{16}$$

has nontrivial solutions. By the definition (9) of $T(\mu)$ this equation is the weak form of the eigenvalue problem

$$\Delta u + \lambda u = 0 \quad \text{in } \Omega = (0,\pi) \times (0,\pi)$$

with the boundary condition (3). To solve this problem with the rule of separating variables, we take the ansatz $u(x,y) = u_1(x)u_2(y) \not\equiv 0$ and derive

$$\frac{u_1''}{u_1} + \frac{u_2''}{u_2} + \lambda = 0 \quad \text{in } \Omega = (0,\pi) \times (0,\pi).$$

Hence, $u_1$ satisfies the equation

$$u_1'' + k_1 u_1 = 0 \quad \text{for some } k_1 \in \mathbf{R} \tag{17}$$

with the boundary conditions $u_1'(0) = u_1'(\pi) = 0$. Similarly, $u_2$ is a solution of

$$u_2'' + k_2 u_2 = 0 \quad \text{for some } k_2 \in \mathbf{R} \tag{18}$$

with the boundary conditions

$$\begin{aligned}
h_0(\mu)u_2(0) - h_1(\mu)u_2'(0) &= 0, \\
h_0(\mu)u_2(\pi) + h_1(\mu)u_2'(\pi) &= 0.
\end{aligned}$$

These are eigenvalue problems of the one-dimensional differential operator $\frac{d^2}{dx^2}$ with two different boundary conditions. Solutions of these problems are of the form

$$\begin{aligned}
(u_1, k_1) &= (\cos nx,\ n), \quad n \in \mathbf{N}, \\
(u_2, k_2) &= \Big(h_0(\mu)\sin(k(\mu)y) + h_1(\mu)k(\mu)\cos(k(\mu)y),\ k(\mu)\Big),
\end{aligned}$$

where $k(\mu) \in \mathbf{R}$ satisfies

$$2h_0(\mu)h_1(\mu)k\cos(k\pi) + (h_0^2(\mu) - h_1^2(\mu)k^2)\sin(k\pi) = 0. \tag{19}$$

Thus eigenvalues of the Laplacian $-\Delta$ are given as

$$\lambda = n^2 + k(\mu)^2 \tag{20}$$

with the corresponding eigenfunction

$$\phi(\mu) := \tilde{\phi}/\|\tilde{\phi}\|, \qquad \tilde{\phi} := (\cos nx)[h_0(\mu)\sin(k(\mu)y) + h_1(\mu)k(\mu)\cos(k(\mu)y)]. \tag{21}$$

As a function of the homotopy parameter $\mu$, the function $\kappa(\mu)$ has the following properties (cf. Mei/Theil [13]).

**Lemma 2** *Under the assumption (4), the equation (19) does not have integer solution for $\mu \in (0,1)$. Moreover, if $[h_1(\mu)/h_0(\mu)]' < 0$ for all $\mu \in (0,1)$, then the solution $k(\mu)$ of (19) increases monotonously from $m \in \mathbf{N}$ at $\mu = 0$ to $(m+1) \in \mathbf{N}$ at $\mu = 1$.*

Note that after multiplying the factor $\sin(k\pi)$ to the both sides of the equation (19) we can rewrite it as

$$[h_0(\mu)\sin(k\pi) - h_1(\mu)k(1-\cos(k\pi))][h_0(\mu)\sin(k\pi) + h_1(\mu)k(1+\cos(k\pi))] = 0.$$

We use the following definition of parities of $k(\mu)$ in Ashwin/Mei [1], which are consistent with the parities of wavenumbers of the Neumann problem at $\mu = 0$.

**Definition 3** *The parities of the wavenumber $k(\mu)$ for Robin boundary conditions are defined as*

$$\kappa(\mu) = \begin{cases} \text{EVEN} & \text{if } h_0(\mu)\sin(k(\mu)\pi) = h_1(\mu)k(\mu)[1 - \cos(k(\mu)\pi)], \\ & \qquad\qquad\qquad\qquad\qquad\qquad\qquad\qquad \text{for all } \mu \in [0,1], \\ \text{ODD} & \text{if } h_0(\mu)\sin(k(\mu)\pi) = -h_1(\mu)k(\mu)[1 + \cos(k(\mu)\pi)]. \end{cases} \tag{22}$$

In the rest of this paper we restrict the discussion to the case

$$\frac{h_0(\mu)}{h_1(\mu)} \geq 0, \qquad \left(\frac{h_1(\mu)}{h_0(\mu)}\right)' < 0.$$

We conclude that bifurcation points of (10) on the trivial solution manifold are

$$\Big\{(0,\ \lambda(\mu),\ \mu);\ \ \lambda(\mu) = n^2 + k(\mu)^2,\ n \in \mathbf{N},\ \text{and}\ k(\mu)\ \text{satisfying (19)}\Big\}, \tag{23}$$

see Figure 2. The kernel $\mathrm{Ker}(D_u G(0, \lambda(\mu), \mu))$ is generically one-dimensional and

$$\mathrm{Ker}(D_u G(0, \lambda(\mu), \mu)) = \mathrm{span}[\phi(\mu)].$$

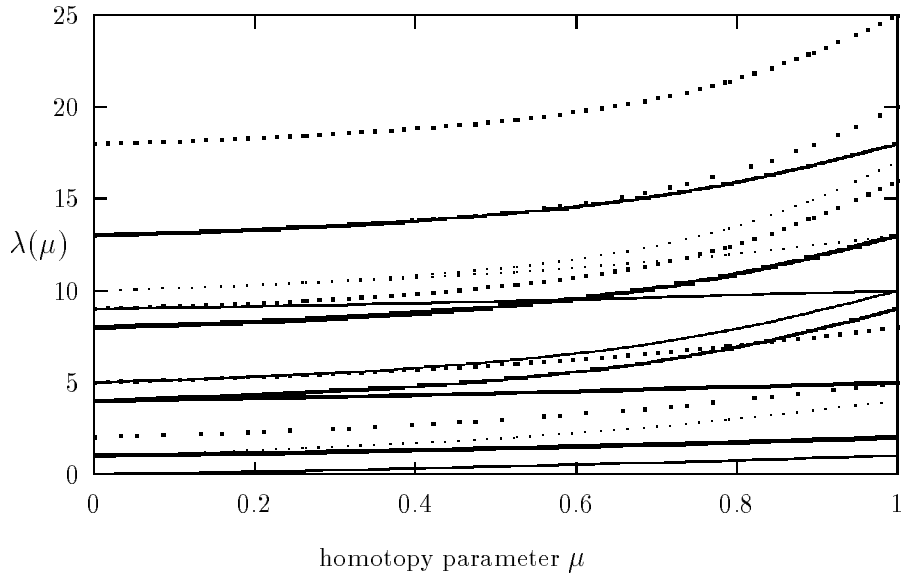

Figure 2: *Bifurcation points of (10) in the parameter space $(\lambda, d)$. Here we have chosen $h_0(\mu) = \mu$, $h_1(\mu) = 1 - \mu$.*

From the statement (20) and the Figure 2 it is evident that two curves of bifurcation points may intersect. In fact, these are generic as $\mu$ approaches zero and one, respectively. A intersection point corresponds to a multiple bifurcation of (10). At $\mu = 0, 1$ solutions of the equation (19) are integers. In particular, at $\mu = 0$ the boundary conditions (3) reduce to homogeneous Neumann type. Thus the eigenvalues of $-\Delta$ are of the form
$$\lambda_0 = n^2 + k^2, \quad n, \ k \in \mathbf{N} \cup \{0\}.$$
Such an eigenvalue is generically double and the associated eigenspace is
$$\mathrm{Ker}(D_u G(0, \lambda_0, 0)) = \mathrm{span}[\phi_1, \phi_2],$$
where
$$\phi_1 = \frac{2}{\pi} \cos(nx) \cos(ky), \quad \phi_2 = \frac{2}{\pi} \cos(kx) \cos(ny).$$
Consequently, as $\mu \to 0$, two different curves of simple bifurcation points approach the same double bifurcation point. However, with $\lambda$ as the bifurcation parameter, there is exactly one solution branch bifurcating from the trivial solution curve at every simple bifurcation point, while at a double bifurcation point for $\mu = 0$ the Neumann problem have up to four different nontrivial solution branches (cf. Mei [11]). On the other hand, as we have seen before, the $D_4$-symmetry of the Neumann problem breaks into $D_2 \times D_2$-symmetry as $\mu$ moves away from 0. This leads to the question how the bifurcation scenario reacts to such a symmetry-breaking in the boundary conditions. We want to investigate which solution branches of the Neumann problem and what symmetries of the solutions persist as the parameter $\mu$ varies.

## 4  Liapunov-Schmidt Reduction

To investigate solutions of (10) at a bifurcation scenario point $(0, \lambda_0, \mu_0)$ on the curve $(0, \lambda(\mu), \mu)$ in (23), we use the well-known Liapunov-Schmidt method to reduce the problem (10) to an algebraic system (cf. Golubitsky/Schaeffer [5] and Golubitsky/Stewart/Schaeffer [6]).

According to the Fredholm properties of $D_u G(0, \lambda_0, \mu_0) = I + (\lambda_0 + 1)T(\mu_0)$, we have the decomposition

$$X = \mathrm{Ker}(D_u G(0, \lambda_0, \mu_0)) \oplus \mathrm{Im}(D_u G(0, \lambda_0, \mu_0)).$$

Suppose that the kernel $\mathrm{Ker}(D_u G(0, \lambda_0, \mu_0)) = \mathrm{span}[\phi_1, \ldots, \phi_\ell]$ is $\ell$-dimensional ($\ell = 1, 2$ generically). We write elements $(u, \lambda, \mu) \in X \times \mathbf{R} \times \mathbf{R}$ as

$$\begin{aligned} u &= \sum_{i=1}^{\ell} z_i \phi_i + w = z \cdot \phi + w, \\ \lambda &= \lambda_0 + \sigma, \\ \mu &= \mu_0 + \nu, \end{aligned}$$

where $z = (z_1, \ldots, z_\ell)$, $\phi = (\phi_1, \ldots, \phi_\ell)$, $z_i$, $\sigma$, $\nu \in \mathbf{R}$ and $w \in \mathrm{Im}(D_u G(0, \lambda_0, \mu_0))$. Consider the projection $Q := I - \sum_{i=1}^{\ell} \langle \phi_i, \cdot \rangle \phi_i$ from $X$ onto $\mathrm{Im}(D_u G(0, \lambda_0, \mu_0))$. We rewrite the equation $G(u, \lambda, \mu) = 0$ into a system

$$QG(z \cdot \phi + w, \ \lambda_0 + \sigma, \ \mu_0 + \nu) = 0, \qquad (24)$$

$$(I - Q)G(z \cdot \phi + w, \ \lambda_0 + \sigma, \ \mu_0 + \nu) = 0. \qquad (25)$$

Solving $w$ uniquely from (24) as a function of $z$, $\sigma$, $\nu$ and substituting it into (25), we obtain the *reduced bifurcation equation* for $z$, $\sigma$, $\mu$

$$(I - Q)G(z \cdot \phi + w(z, \sigma, \nu), \ \lambda_0 + \sigma, \ \mu_0 + \nu) = 0. \qquad (26)$$

We note that $w(0, 0, 0) = 0$, $D_z w(0, 0, 0) = 0$ from (24). In the coordinate system

$$z \cdot \phi \in \mathrm{Ker}(D_u G(0, \lambda_0, \mu_0)) \longleftrightarrow z = (z_1, \ldots, z_\ell) \in \mathbf{R}^\ell,$$

the operator equation (26) becomes a system of $\ell$ algebraic equations

$$B(z, \sigma, \nu) := \Big(\langle \phi_i, \ G(z \cdot \phi + w(z, \sigma, \nu), \ \lambda_0 + \sigma, \ \mu_0 + \nu)\rangle\Big)_{i=1}^{\ell} = 0. \tag{27}$$

By definition the projection $Q$ is $\Gamma$-equivariant. Thereafter the mapping $B(z, \sigma, \nu)$ is also $\Gamma$-equivariant with respect to the induced action of $\Gamma$ in $\mathbf{R}^\ell$, i.e.,

$$B(\gamma z, \sigma, \nu) = \gamma B(z, \sigma, \nu) \quad \text{for all} \ \ (z, \sigma, \nu) \in \mathbf{R}^\ell \times \mathbf{R} \times \mathbf{R}, \ \gamma \in \Gamma.$$

We take the Taylor expansion of the components of $(\lambda_0 + 1)B(z, \sigma, \nu) = 0$ at the point $(z, \sigma, \nu) = (0, 0, 0)$ and consider the truncated form

$$-\sigma z_i + (\lambda_0 + 1)^2 \langle \phi_i, \ T'(\mu_0)(z \cdot \phi)\rangle \nu \tag{28}$$
$$+ \ \Big\langle \phi_i, \ \frac{1}{2}D_{uu}f_0(z \cdot \phi)^2 + D_{uu}f_0(z \cdot \phi)\Big(\frac{1}{2}\sum_{|\alpha|=2}D^\alpha w_0 z^\alpha\Big) + \frac{1}{6}D_{uuu}f_0(z \cdot \phi)^3\Big\rangle$$
$$= \ 0, \ \ i = 1, \ldots, \ell.$$

Here and in the sequel $D_{uu}f_0$ and $D_{uuu}f_0$ denote the derivatives of $f$ at $(u, \lambda) = (0, \lambda_0)$; $\alpha \in \mathbf{N}^\ell$ is a multi-index and

$$\frac{1}{2}\sum_{|\alpha|=2}D^\alpha w_0 z^\alpha$$

represents the second order terms in the Taylor expansion of $w$ at $(z, \sigma, \nu) = (0, 0, 0)$.

We recall that in the singularity theory a problem $B = 0$ is finitely determined if there exists $k \in \mathbf{N}$, $k < \infty$, such that the bifurcation scenario of $B = 0$ and its $k$-jets $j_k(B) = 0$, the Taylor expansion of $B$ truncated at $k$-th order, are equivalent. The determinacy of a general problem at a bifurcation point is characterized by its reduced bifurcation equations. We refer to Golubitsky/Schaeffer [5] for more detailed discussions. For 3-determined problems solutions of the system (28) correspond one-to-one to those of the original problem (10) and contain all information of bifurcations of (10) at $(0, \lambda_0, \mu_0)$.

With the knowledge of the bifurcation point $(0, \lambda_0, \mu_0)$ and the kernel $\text{Ker}(D_u G_0)$, the terms $\langle \phi_i, \ D_{uu}f_0(z \cdot \phi)^2\rangle$ and $\langle \phi_i, D_{uuu}f_0(z \cdot \phi)^3\rangle$ in (28) can be calculated directly. The other terms involve the derivatives $T'(\mu_0)$ and $D_{z_i z_j}w_0$. Since the function $w(z, \sigma, \nu)$ is defined implicitly by the equation (24), the term $D_{z_i z_j}w_0$ is described as the unique solution $v$ of the linear problem

$$D_u G_0 v = QT(\mu_0)D_{uu}f_0\phi_i\phi_j, \quad v \in \text{Im}(D_u G_0). \tag{29}$$

The term $T'(\mu_0)(z \cdot \phi)$ is calculated as a solution of the equation (11). In fact, via (13) we derive

$$\begin{aligned}
& \langle \phi_i, \ T'(\mu_0)(z \cdot \phi) \rangle \\
=& \ \tilde{h}(\mu_0) \langle \phi_i, \ \left(-\frac{y^2}{\pi} + y\right) T(\mu_0)(z \cdot \phi) \\
& + T(\mu_0)\left[\frac{2}{\pi} T(\mu_0)(z \cdot \phi) + 2\left(\frac{2y}{\pi} - 1\right)\frac{\partial}{\partial y}\left(T(\mu_0)(z \cdot \phi)\right) + \left(\frac{y^2}{\pi} - y\right)(z \cdot \phi)\right]\rangle \\
=& \ \tilde{h}(\mu_0) \frac{1}{(\lambda_0+1)^2} \left\langle \phi_i, \ \frac{2}{\pi}(z \cdot \phi) - 2(\lambda_0+1)\left(\frac{2y}{\pi}-1\right)\frac{\partial}{\partial y}\left(T(\mu_0)(z \cdot \phi)\right)\right\rangle \\
=& \ \tilde{h}(\mu_0) \frac{2}{(\lambda_0+1)^2}\left[\frac{z_i}{\pi} + \left\langle \phi_i, \ \left(\frac{2y}{\pi}-1\right)\frac{\partial}{\partial y}(z \cdot \phi)\right\rangle\right],
\end{aligned}$$

where $\tilde{h}(\mu_0) = \left(\dfrac{h_0}{h_1}\right)'(\mu_0)$.

# 5 Bifurcation Scenarios

## 5.1 Simple bifurcations for $\mu \in (0, 1)$

Let $\lambda(\mu) = n^2 + k(\mu)^2$ be a homotopy of simple eigenvalues of the Laplacian and $\lambda(0) = n^2 + m^2$, $\lambda(1) = n^2 + (m+1)^2$. Then $(0, \lambda(\mu), \mu)$ is a curve of simple bifurcation points of (10) and

$$\mathrm{Ker}(D_u G)(0, \lambda(\mu), \mu) = \mathrm{span}[\phi]$$

is 1-dimensional and $\phi$ is given in (21). To obtain the generic bifurcation diagram of (10) at $(0, \lambda(\mu_0), \mu_0)$ for an arbitrary $\mu_0 \in (0, 1)$, we consider the equation (28), i.e., 3-jet of the reduced bifurcation equation,

$$\begin{aligned}
0 =& \ -\sigma z + (\lambda_0+1)^2 \langle \phi, \ T'(\mu_0)\phi\rangle \nu z \\
& + \left\langle \phi, \ \frac{1}{2} D_{uu} f_0 \phi^2 \right\rangle z^2 + \left\langle \phi, \ \frac{1}{2} D_{uu} f_0 (D_{zz} w_0)\phi + \frac{1}{6}(D_{uuu} f_0)\phi^3 \right\rangle z^3.
\end{aligned} \tag{30}$$

Here $z \in \mathbf{R}$ is a scalar and

$$\langle \phi, \ T'(\mu_0)\phi \rangle = \left(\frac{g_0}{g_1}\right)'(\mu_0) \frac{2}{(\lambda_0+1)^2}\left[\frac{1}{\pi} + \left\langle \phi, \ \left(\frac{2y}{\pi}-1\right)\frac{\partial}{\partial y}\phi\right\rangle\right].$$

**Theorem 4** *The problem (10) undergoes a pitchfork bifurcation at all points on the curve $(0, \lambda(\mu), \mu)$, $\mu \in (0, 1)$, i.e, the simple bifurcation points. Moreover, the truncated bifurcation equation (30) reduces to*

$$j_3[(\lambda_0+1)B(z, \sigma, \nu)] = (-\sigma + a\nu)z + cz^3 = 0, \tag{31}$$

*with*

$$a = 2\tilde{h}(\mu_0)\left[\frac{1}{\pi} + \left\langle \phi, \left(\frac{2y}{\pi} - 1\right)\frac{\partial}{\partial y}\phi\right\rangle\right],$$
$$c = \left\langle \phi, \frac{1}{2}D_{uu}f_0\left(D_{zz}w_0\right)\phi + \frac{1}{6}D_{uuu}f_0\phi^3\right\rangle.$$
(32)

**Proof:** It is easy to verify that the eigenfunction $\phi$ of the Laplacian has the property $\langle \phi, \phi^2 \rangle = 0$. Thus the second order term in (30) vanishes. The conclusion follows directly from equations (30) and (31) consecutively. ∎

The nontrivial solution of (31) is given as

$$z = \left(\frac{\sigma - a\nu}{c}\right)^{1/2}.$$

## 5.2 Double bifurcations of the Neumann problem

For the Neumann problem ($\mu = 0$) a generic double bifurcation point $(0, \lambda_0, 0)$ has the property $\lambda_0 = n^2 + k^2(0)$ with the wavenumbers $n$, $k := k(0)$ ($n \neq k$) as integers. Furthermore, we can choose

$$\text{Ker}(D_u G_0) = \text{span}[\phi_1, \phi_2]$$

with

$$\phi_1 := \begin{cases} \dfrac{2}{\pi}\cos(nx)\cos(ky), & \text{for } n \cdot k \neq 0; \\ \dfrac{\sqrt{2}}{\pi}\cos(nx), & \text{for } n \neq 0, \ k = 0; \\ \dfrac{\sqrt{2}}{\pi}\cos(ky), & \text{for } n = 0, \ k \neq 0, \end{cases}$$

and

$$\phi_2 := \phi_1(y, x).$$

On the other hand, taking into account the homotopy parameter $\mu$ in the boundary conditions, we see this double bifurcation point is split into two simple bifurcation points $(0, \lambda_i(\mu), \mu)$, $i = 1, 2$ with $\lambda_1 = n^2 + k^2(\mu)$ and $\lambda_2 = k^2 + n^2(\mu)$ for $\mu \neq 0$. We want to investigate bifurcation scenario of (10) at a double bifurcation point $(0, \lambda_0, 0)$ and its variation with respect to the homotopy parameter $\mu$.

Note that $\langle \phi_i, \phi_j \phi_l \rangle = 0$ for all $i, j, l = 1, 2$. The equation (29) can be solved analytically (cf. Mei [11]). Together with the statements

$$\langle \phi_i^4, 1 \rangle = \begin{cases} \dfrac{9}{4\pi^2} & \text{for } n \cdot k \neq 0, \\ \dfrac{3}{2\pi^2} & \text{for } n \cdot k = 0, \ n^2 + k^2 \neq 0 \end{cases}$$

and $\langle \phi_i^2, \phi_j^2 \rangle = \dfrac{1}{\pi^2}$ for $i \neq j$, we simplify the equations (28) into

$$\begin{aligned}
-\sigma z_1 + (\lambda_0 + 1)^2 \langle \phi_1, T'(\mu_0)(z_1\phi_1 + z_2\phi_2)\rangle \nu + c_1 z_1^3 + c_2 z_1 z_2^2 &= 0, \\
-\sigma z_2 + (\lambda_0 + 1)^2 \langle \phi_2, T'(\mu_0)(z_1\phi_1 + z_2\phi_2)\rangle \nu + c_2 z_1^2 z_2 + c_1 z_2^3 &= 0.
\end{aligned} \quad (33)$$

Here $c_1$, $c_2$ are constants. More precisely, if $n \cdot k \neq 0$, we have

$$\begin{aligned}
c_1 &= \frac{1}{6\pi^2}\left[\frac{9}{4}D_{uuu}f_0 - \frac{1}{4}(D_{uu}f_0)^2 \frac{45(k^2-n^2)^2 + 4k^2n^2}{(k^2-3n^2)(n^2-3k^2)(n^2+k^2)}\right], \\
c_2 &= \frac{1}{6\pi^2}\left[3D_{uuu}f_0 - 6(D_{uu}f_0)^2 \frac{1}{n^2+k^2}\left(\frac{(k^2-n^2)^2 - 4k^2n^2}{[(k^2+n^2)^2 - 16k^2n^2]} - \frac{1}{2}\right)\right].
\end{aligned}$$

If $n = 0$, $k \neq 0$, then

$$c_1 = \frac{1}{6\pi^2}\left[\frac{3}{2}D_{uuu}f_0 + \frac{5}{2k^2}(D_{uu}f_0)^2\right],$$

$$c_2 = \frac{1}{2\pi^2}D_{uuu}f_0.$$

If $\nu = 0$, the equations (33) coincide with those in Mei [11], and yield four nontrivial solutions of (10) with symmetries of the isotropy subgroups of $\phi_1$, $\phi_2$ and $\phi_1 \pm \phi_2$, respectively. For $\nu \neq 0$, the forced symmetry-breaking in boundary conditions introduce in (33) the terms

$$\begin{aligned}
& (\lambda_0 + 1)^2 \langle \phi_i, T'(\mu_0)(z_1\phi_1 + z_2\phi_2)\rangle \nu \\
&= 2\left(\frac{g_0}{g_1}\right)'(0)\left[\frac{z_i}{\pi} + \langle \phi_i, (\frac{2y}{\pi} - 1)\frac{\partial}{\partial y}(z_1\phi_1 + z_2\phi_2)\rangle\right] \\
&=: d_i z_i, \quad i = 1, 2.
\end{aligned}$$

Here,

$$d_1 = d_2 = \frac{4}{\pi}\left(\frac{g_0}{g_1}\right)'(0), \quad \text{if } n \cdot k \neq 0,$$

$$d_1 = \frac{4}{\pi}\left(\frac{g_0}{g_1}\right)'(0), \quad d_2 = 0, \quad \text{if } n = 0, \, k \neq 0.$$

The system (33) reduces to

$$\begin{aligned}
\left[-\sigma + d_1\nu + c_1 z_1^2 + c_2 z_2^2\right]z_1 &= 0, \\
\left[-\sigma + d_2\nu + c_2 z_1^2 + c_1 z_2^2\right]z_2 &= 0.
\end{aligned} \quad (34)$$

**Remark:** The coefficients in the equations (32) and (34) are related as follows have properties

$$\lim_{\mu \to 0} a = d_i \quad \text{for } \phi = \phi_i, \ i = 1, 2$$
$$\lim_{\mu \to 0} c = c_1.$$

Solutions of the system (34) are

$$\begin{aligned}
a) & \quad \left(\pm \left(\frac{\sigma - d_1 \nu}{c_1}\right)^{1/2}, \ 0\right), \quad \left(0, \ \pm \left(\frac{\sigma - d_2 \nu}{c_1}\right)^{1/2}\right) \\
b) & \quad \left(\pm \left(\frac{(c_1 - c_2)\sigma - (c_1 d_1 - c_2 d_2)\nu}{c_1^2 - c_2^2}\right)^{1/2}, \ \pm \left(\frac{(c_1 - c_2)\sigma + (c_2 d_1 - c_1 d_2)\nu}{c_1^2 - c_2^2}\right)^{1/2}\right).
\end{aligned} \quad (35)$$

These lead to four bifurcating solution branches of the original problem (10), i.e., (1). The solutions in (35a) are pure mode solution branches with the isotropy groups of $\phi_1$, $\phi_2$, respectively. They correspond to those bifurcating at the simple bifurcation points on the curves $(0, \lambda_i(\mu), \mu)$, $i = 1, 2$. The solutions in (35b) involve both $\phi_1$ and $\phi_2$ modes. They are called the mixed mode branches.

If $n \cdot k = 0$, the terms $d_1 \nu$, $d_2 \nu$ break the $D_4$-symmetry of the Neumann problem and the mixed mode solution branches have merely the trivial symmetry. Moreover, the pure mode and mixed mode solution branches may intersect at

$$\sigma = \frac{(c_1 d_1 - c_2 d_2)\nu}{(c_1 - c_2)} \quad \text{or} \quad \sigma = \frac{(c_2 d_1 - c_1 d_2)\nu}{(c_1 - c_2)},$$

and induce a secondary bifurcation, respectively.

If $n \cdot k \neq 0$, then $d_1 = d_2$ and the $D_4$-symmetry is preserved in (34) and the mixed mode solutions (35b) becomes

$$\pm \left(\left(\frac{\sigma - d_1 \nu}{c_1 + c_2}\right)^{1/2}, \ \pm \left(\frac{\sigma - d_1 \nu}{c_1 + c_2}\right)^{1/2}\right).$$

Symmetries of these solutions are the isotropy groups of the eigenfunctions $\phi_1 + \phi_2$ and $\phi_1 - \phi_2$, respectively.

As an conclusion, we see that all four bifurcating solution branches of the Neumann problem at a double bifurcation point persist if we vary both $\lambda$ and $\mu$ as bifurcation parameters. Moreover, symmetry of these bifurcating solution branches is preserved for those with the wavenumbers $n \cdot k \neq 0$ and is broken for those with the wavenumbers $n = 0$ or $k = 0$.

# 6  A Simple Example

Choose
$$f(u, \lambda) = \lambda(u^2 + u^3). \tag{36}$$

We consider bifurcation scenarios at the corank-2 bifurcation points $u_0 = 0$, $\mu_0 = 0$ and $\lambda_0 = 5$, 10, 20, respectively. Moreover, we examine variations of the bifurcation scenarios as homotopy parameter $\mu$ moves away from zero, i.e., as the homogeneous Neumann boundary conditions with $D_4$-symmetry are perturbed. To simplify the discussion, we take $h_0(\mu) = \mu$, $h_1(\mu) = 1 - \mu$. Note that at $\mu = 0$ we have

$$\mathrm{Ker}(D_u G_0) = \mathrm{span}[\phi_1, \phi_2].$$

and the inequalities $c_1 \neq 0$, $c_2 \neq 0$ and $c_1^2 - c_2^2 \neq 0$ holds for all $n$, $k \in \mathbf{N} \cup \{0\}$.

**(1) Wavenumbers** $n = 1$, $k = 2$: For $\lambda(0) = 5$ and $n = 1$, $k = 2$ we have

$$\phi_1 = \frac{2}{\pi} \cos x \cos 2y, \qquad \phi_2 = \frac{2}{\pi} \cos 2x \cos y.$$

Furthermore, $d = \dfrac{4}{\pi}$, $c_1 = \dfrac{5695}{132\pi^2}$, $c_2 = \dfrac{110220}{132\pi^2}$. Solutions $(z_1, z_2)$ in (35) becomes

$$\left(\pm \Big(\frac{132\pi(\pi\sigma - 4\nu)}{5695}\Big)^{1/2},\ 0\right), \qquad \left(0,\ \pm \Big(\frac{132\pi(\pi\sigma - 4\nu)}{5695}\Big)^{1/2}\right)$$

$$\pm \left(\Big(\frac{132\pi(\pi\sigma - 4\nu)}{115915}\Big)^{1/2},\ \pm \Big(\frac{132\pi(\pi\sigma - 4\nu)}{115915}\Big)^{1/2}\right).$$

Figures 3 and 4 show the pure and mixed modes solution branches.

**(2) Wavenumbers** $n = 0$, $k = 1$: We have $\lambda(0) = 1$ for $n = 0$, $k = 1$, moreover,

$$\phi_1 = \frac{\sqrt{2}}{\pi} \cos y, \qquad \phi_2 = \frac{\sqrt{2}}{\pi} \cos x.$$

Simple calculations show $d_1 = \dfrac{4}{\pi}$, $d_2 = 0$ $c_1 = \dfrac{19}{6\pi^2}$, and $c_2 = \dfrac{3}{\pi^2}$. The solution branches described by $(z_1, z_2)$ in (35) are

$$\left(\pm \Big(\frac{6\pi(\pi\sigma - 4\nu)}{19}\Big)^{1/2},\ 0\right), \qquad \left(0,\ \pm \Big(\frac{6\pi^2\sigma}{19}\Big)^{1/2}\right),$$

$$\pm \left(\Big(\frac{6\pi(\pi\sigma - 76\nu)}{37}\Big)^{1/2},\ \pm \Big(\frac{6\pi(\pi\sigma + 72\nu)}{37}\Big)^{1/2}\right).$$

The pure $\phi_2$-mode solution branch meets a mixed mode solution at $\sigma = 76\nu/\pi$ and induces a secondary bifurcation.

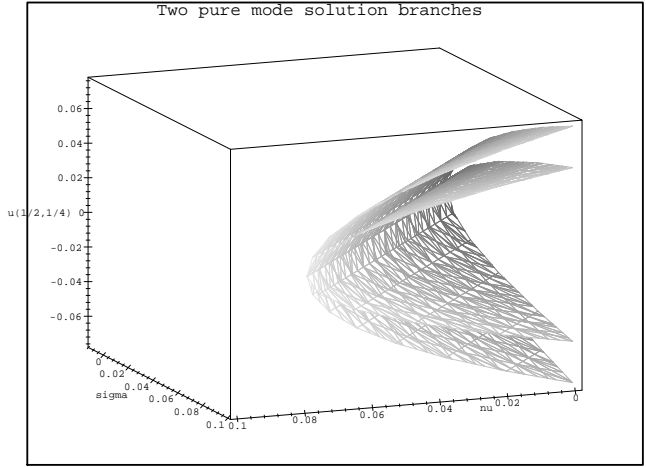

Figure 3: *Two solution branches of pure $\phi_1$ and $\phi_2$ modes.*

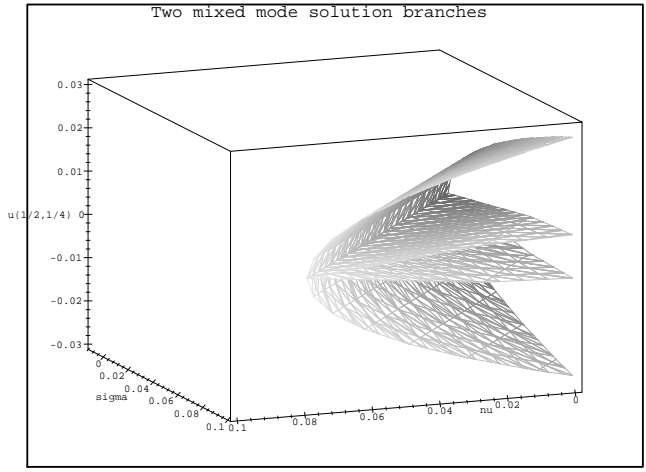

Figure 4: *Two solution branches of mixed $\phi_1$ and $\phi_2$ modes.*